# THE MIXING TIME FOR SIMPLE EXCLUSION

By Ben Morris

*University of California, Davis*

We obtain a tight bound of $O(L^2 \log k)$ for the mixing time of the exclusion process in $\mathbf{Z}^d/L\mathbf{Z}^d$ with $k \leq \frac{1}{2} L^d$ particles. Previously the best bound, based on the log Sobolev constant determined by Yau, was not tight for small $k$. When dependence on the dimension $d$ is considered, our bounds are an improvement for all $k$. We also get bounds for the relaxation time that are lower order in $d$ than previous estimates: our bound of $O(L^2 \log d)$ improves on the earlier bound $O(L^2 d)$ obtained by Quastel. Our proof is based on an auxiliary Markov chain we call the *chameleon process*, which may be of independent interest.

**1. Introduction.** Let $G = (V, E)$ be a finite, connected graph and define a *configuration* as follows. In a configuration, each vertex in $V$ contains either a black ball or a white ball (where balls of the same color are indistinguishable), and the number of black balls is at most $|V|/2$. The *exclusion process* on $G$ is the following continuous-time Markov process on configurations. For each edge $e$ at rate 1, switch the balls at the endpoints of $e$.

Note that since the exclusion process is irreducible and has symmetric transition rates, the uniform distribution is stationary. Let $\mathcal{C}$ denote the space of configurations and for probability distributions $\mu, \nu$ on $\mathcal{C}$, let

$$\|\mu - \nu\| = \max_{Q \subset \mathcal{C}} \mu(Q) - \nu(Q) = \min_{X \sim \mu, Y \sim \nu} \mathbf{P}(X \neq Y)$$

be the total variation distance. Denote the mixing time by

(1) $\qquad \tau_{\text{mix}} = \inf\{t : \|K^t(x, \cdot) - \mathcal{U}\| \leq \tfrac{1}{4} \text{ for all } x \in \mathcal{C}\},$

where $K$ is the transition kernel for the exclusion process and $\mathcal{U}$ is the uniform distribution over configurations.

This paper is concerned with bounding the mixing time in the important special case where $G$ is the $d$-dimensional torus $\mathbf{Z}^d/L\mathbf{Z}^d$; we call this process *simple exclusion*.

The exclusion process is a widely studied Markov chain, with connections to card shuffling [1, 14], statistical mechanics [2, 6, 11, 13] and a variety of other processes (see, e.g., [5, 8]); it has been one of the major examples driving the study of mixing times for Markov chains (see, e.g., [1, 3, 5, 14]).

Quastel [11] and a subsequent independent paper by Diaconis and Saloff-Coste [3] bounded the spectral gap for simple exclusion by comparing the process with Bernoulli–Laplace diffusion (i.e., the exclusion process on the complete graph), obtaining a bound of $O(L^2 d)$ for the *relaxation time* (i.e., inverse spectral gap); by a standard theorem (see, e.g., [12]), this implies

$$\tau_{\mathrm{mix}} \le c[L^2 d] \log \binom{L^d}{k},$$

where $k$ is the number of black balls [so that $\binom{L^d}{k}$ is the number of states in the chain] and $c$ is a universal constant.

An improved dependence on the number of states can be obtained using the *log Sobolev constant*. In a famous result of Yau [15], which built on the earlier work of Lu and Yau [9], the log Sobolev constant for simple exclusion was shown to be $O(c_d L^2)$, where $c_d$ is a constant that depends only on $d$; by a standard theorem (see, e.g., [12]), this implies

$$(2) \qquad \tau_{\mathrm{mix}} \le c_d L^2 \log \log \binom{L^d}{k}.$$

(An analysis that combines the comparison technique [3, 4] with Lee and Yau's [7] estimate for the log Sobolev constant for Bernoulli–Laplace diffusion shows that we can take $c_d = cd$ for a universal constant $c$ when $k = L^d/2$.)

The bound (2) was the best known, but it was not tight for small $k$; in this case, the log Sobolev constant, which is based on a single analytic inequality, does not capture the faster mixing of the process in its early stages.

The main contribution of this paper is to give a bound of the correct order, namely $O(L^2 \log k)$, for the mixing time. We accomplish this using an auxiliary Markov chain called the *chameleon process*, a variant of the exclusion process in which the balls can change color. Using a relationship between the transition probabilities of the two processes (see Lemma 1 below), we give a detailed analysis of the mixing during the various stages and obtain a tight bound.

As an added benefit, when dependence on the dimension $d$ is considered, our bounds are an improvement for all values of $k$, both for the mixing time and for the relaxation time. Table 1 shows the improvement when $k = L^d/2$.

We have included the special case $G = \{0,1\}^d$ above to illustrate the dependence on $d$. In this case, Wilson [14] showed that the mixing time is at least of order $d$ (see [14]), so our bound is within a factor $O(\log d)$ of the true answer. Diaconis and Saloff-Coste [3] wondered if the factor $d$ in the estimate $\tau_{\text{relax}} \leq cd$ was necessary. We show that it can be replaced by $\log d$.

We note that the analyses of Quastel, Diaconis and Saloff-Coste and Lee and Yau apply to the box $\{0, 1, \ldots, L\}^d$ in addition to the torus. Our results also hold for this and other variants, but in the present paper we will only treat the case where $G$ is the torus.

There is a connection between our key tool—the chameleon process—and *evolving sets*, which were used by the author and Peres [10] to prove heat kernel bounds based on isoperimetric inequalities. The reader who is familiar with [10] might recognize the chameleon process as a variant of an evolving set process (in somewhat disguised form). This connection will not be used explicitly.

## 2. Interchange and chameleon processes.

2.1. *Interchange process.* We will now describe a related Markov chain, the *interchange process*, which will be useful. Let $n = |V|$. In a configuration of the interchange process, balls labeled $1, 2, \ldots, n$ occupy the vertices of $G$. Each vertex contains exactly one ball, which is colored either black, white, red or pink. The transition rule is the following: For each edge $e$ at rate 2, flip a coin. If the the coin lands heads, switch the balls at the endpoints of $e$; else do nothing. More precisely, there are independent Poisson processes of rate 2 for every edge; when the Poisson process corresponding to edge $e$ has an arrival, the balls at the endpoints of $e$ are switched w.p. $1/2$. When there is an arrival that corresponds to $e$, we will say *the clock of edge $e$ rings* or simply *edge $e$ rings*.

For a ball $i \in \{1, \ldots, n\}$, define the *trajectory of $i$* as the map $t \to v_t(i)$, where $v_t(i)$ is the position of ball $i$ at time $t$. Note that since each pair of neighboring balls is switched at rate 1, the *unordered set* $\{v_t(1), \ldots, v_t(k)\}$ behaves like the set of black balls in the exclusion process.

Table 1

| | $\tau_{\text{mix}}$ | | $\tau_{\text{relax}}$ | |
| --- | --- | --- | --- | --- |
| $G$ | Old bound | New bound | Old bound | New bound |
| $\mathbf{Z}^d/L\mathbf{Z}^d$ | $cd^2L^2\log(L)$ | $cd\log(d)L^2\log(L)$ | $cdL^2$ | $c\log(d)L^2$ |
| $\{0,1\}^d$ | $cd^2$ | $cd\log(d)$ | $cd$ | $c\log(d)$ |

2.2. *Chameleon process.* The chameleon process is a variant of the interchange process in which the balls can change color. In the chameleon process, the balls move in the same way as in the interchange process; however, in addition, each time the clock of an edge $e$ rings a two-step *recoloring* operation is performed. The first step, which we call *pinkening*, takes place when the balls at the endpoints of $e$ are red and white, respectively; in this case they are both recolored pink. The second step, which we call *depinking*, takes place when there are a large number of pink balls; in this case all of the pink balls are collectively recolored either red or white with probability 1/2 each. More precisely, the transition rule is the following: When the clock of edge $e$ rings, switch the balls at the endpoints of $e$ if and only if the interchange process switches them. In addition, perform the following two-step recoloring operation:

1. (Pinkening step) If the endpoints of $e$ contain a red and a white ball, respectively, then recolor both balls pink.
2. (Depinking step) If the number of pink balls (including the pink balls created in step 1) is at least as large as *either*

   - the number of red balls or
   - the number of white balls,

   then flip a coin. If the coin lands heads, recolor all pink balls red; else recolor all pink balls white.

*Starting colors.* We will always start with the balls colored in the following way. For some $b \in \{0, \ldots, k-1\}$, balls $1, \ldots, b$ are colored black, ball $b+1$ is colored red and the remaining balls are colored white. Note that the number of black balls $b$ remains fixed for the duration of the process, and that the process is symmetric with respect to the roles of the red and white balls.

2.3. *Overview of the chameleon process.* We now introduce some notation. For a configuration $M$ of the chameleon process, let $R(M)$, $W(M)$ and $P(M)$ denote the set of balls in $M$ that are red, white and pink, respectively, and let

$$r(M) = |R(M)|; \qquad w(M) = |W(M)|; \qquad p(M) = |P(M)|.$$

For a chameleon process $\{M_t, t \geq 0\}$, define $R_t := R(M_t)$ and so on.

We will call the times at which the chameleon process recolors its pink balls *depinking times*. Let $T_0 = 0$ and for $j \in \{1, 2, \ldots\}$, let $T_j$ be the $j$th depinking time. We will now describe the behavior of the chameleon process between time $T_j$ and $T_{j+1}$. At time $T_j$, there are only black, white and red balls (since depinking changes all pink balls to either red or white). Let $b$ be the number of black balls and let $m = n - b$. Assume that there are at least as many white balls as red balls at time $T_j$, so that $r_{T_j} \leq \frac{m}{2}$ (the behavior can be described symmetrically when $r_{T_j} > \frac{m}{2}$). From time $T_j$ until time $T_{j+1}$, the process behaves like the interchange process, except that when an edge connecting a red ball to a white ball rings, the two balls are each recolored pink (in the pinkening step). Pink balls accumulate until there are at least as many pink balls as red balls, and then the next depinking occurs at time $T_{j+1}$ (at which point the pink balls are recolored either red or white). Since each pinkening between time $T_j$ and time $T_{j+1}$ increases the number of pink balls by two and decreases the number of red balls by one, the depinking occurs when the number of red balls falls below two-thirds of its value at time $T_j$. Hence we have

$$r_{T_{j+1}} = r_{T_j} \pm \left\lceil \frac{r_{T_j}}{3} \right\rceil.$$

When there are more red balls than white balls at time $T_j$, the process behaves as described above, with the roles of red and white balls reversed. Hence, in general we have

$$r_{T_{j+1}} = r_{T_j} \pm \Delta(r_{T_j}),$$

where $\Delta(x) := \lceil \frac{x \wedge (m-x)}{3} \rceil$. The decision whether to choose $+$ or $-$ is made by a coin flip that is independent of the trajectories of the balls. Hence, the process $\{r_{T_j} : j \geq 0\}$ is a Markov chain. (Note that it is also a martingale.)

**3. Main theorem.** Define the redness $\rho(u, M)$ of vertex $u$ in configuration $M$ by

$$\rho(u, M) := \begin{cases} 1, & \text{if the ball at } u \text{ is red;} \\ 1/2, & \text{if the ball at } u \text{ is pink;} \\ 0, & \text{otherwise.} \end{cases}$$

If $\{M_t, t \geq 0\}$ is a chameleon process, define $\rho_t(u) := \rho(u, M_t)$. The quantity $\sum_u \rho_t(u) = r_t + \frac{1}{2} p_t$ is a measure of the amount of "red paint" in the system at time $t$, if we think of the pink balls as being painted with a mixture of red and white paint. Similarly, the quantity $w_t + \frac{1}{2} p_t$ measures the amount of white paint. The Markov chain $\{r_{T_j} : j = 0, 1, \ldots\}$ has absorbing states at $0$ and $m$. Let $A$ be the event that it absorbs at $m$, the system eventually consisting of only red and black balls with all "white paint" gone. Since $\{r_{T_j}\}$ is a martingale and $r_{T_0} = r_0 = 1$, the optional stopping theorem implies $P(A) = \frac{1}{m}$. Note that $A$ is completely determined by the coin flips that are performed during the depinking steps and hence is independent of the trajectories of the balls.

Define $v_t(1, \ldots, k) := (v_t(1), \ldots, v_t(k))$ and let $X_t = v_t(1, \ldots, k)$ be the vector of positions of balls $1, \ldots, k$ at time $t$. Since the unordered set $\{v_t(1), \ldots,$ $v_t(k)\}$ behaves like the set of black balls in the exclusion process, any upper bound on the mixing time of the process $\{X_t, t \geq 0\}$ also applies to the $k$-particle exclusion process.

Let $\mu_t$ denote the distribution of $X_t$. Then $\mu_t \to \mathcal{U}$ as $t \to \infty$, where $\mathcal{U}$ is the distribution of $k$ uniform samples without replacement from $V$.

For $j \leq k$, we write

$$
\begin{aligned}
\mu_t(u_1,\ldots,u_j) &:= \mathbf{P}(v_t(1,\ldots,j) = (u_1,\ldots,u_j)); \\
\mu_t(u_j|u_1,\ldots,u_{j-1}) &:= \mathbf{P}(v_t(j) = u_j | v_t(1,\ldots,j-1) = (u_1,\ldots,u_{j-1})),
\end{aligned}
\tag{3}
$$

with analogous notation for $\mathcal{U}$. Note that each of the conditional distributions $\mu_t(\cdot|u_1,\ldots,u_{j-1})$ converges to uniform as $t \to \infty$. When all of these distributions are "very close" to uniform, then $\mu_t$ is close to $\mathcal{U}$. (In fact, we only need the conditional distributions to be close, "on average"; Lemma 12 in Section 5.3 formalizes this.) Lemma 2 below relates the mixing of each of the conditional distributions $\mu_t(\cdot|u_1,\ldots,u_b)$ to a corresponding chameleon process with $b$ black balls. It gives a bound based on the expected amount of "white paint" in the system at time $t$, conditional on $A$.

To prove Lemma 2, we will need the following lemma, which indicates the fundamental relationship between the interchange process and the chameleon process. Recall that $\rho_t(u)$ denotes the redness of vertex $u$ at time $t$.

LEMMA 1. *Fix $b \in \{0,\ldots,k-1\}$ and let $z_1,\ldots,z_{b+1}$ be a sequence of distinct vertices in $V$. Consider a chameleon process with $b$ black balls. Then*

$$\mu_t(z_{b+1}|z_1,\ldots,z_b) = \mathbf{E}(\rho_t(z_{b+1})|v_t(1,\ldots,b) = (z_1,\ldots,z_b)). \tag{4}$$

PROOF. Let $\tau_1 < \tau_2 < \cdots$ be the times when edges ring. For $j \geq 1$, let $\{M_t^j, t \geq 0\}$ be the process that behaves like the chameleon process before time $\tau_j$ and behaves like the interchange process at time $\tau_j$ and afterward (i.e., all recoloring is suppressed from time $\tau_j$ onward). Note that $\{M_t^1 : t \geq 0\}$ is the interchange process. We claim that for every $j \geq 1$ we have

$$\mu_t(z_{b+1}|z_1,\ldots,z_b) = \mathbf{E}(\rho(z_{b+1}, M_t^j)|v_t(1,\ldots,b) = (z_1,\ldots,z_b)), \tag{5}$$

where $\rho(z_{b+1}, M_t^j)$ denotes the redness of vertex $z_{b+1}$ in configuration $M_t^j$. This implies (4) because $M_t$ is the a.s. limit of $M_t^j$ as $j \to \infty$. We will show that (5) holds for all $j \geq 1$ by induction. The base case $j = 1$ is trivial because $\{M_t^1\}$ is the interchange process and has no recoloring, so ball $b+1$ is always the only red ball (and there are never any pink balls). Thus, it remains to show that for every $j \geq 1$ we have

$$\mathbf{E}(\rho(z_{b+1}, M_t^j)|D) = \mathbf{E}(\rho(z_{b+1}, M_t^{j+1})|D), \tag{6}$$

where $D$ is the event that $v_t(1,\ldots,b) = (z_1,\ldots,z_b)$. To show this, we introduce a new process $\{M_t^{j+}\}$ that behaves like $\{M_t^j\}$ except that it performs a pinkening step at time $\tau_{j+1}$ (but no depinking then). We first show that

(7) $$\mathbf{E}(\rho(z_{b+1}, M_t^j)|D) = \mathbf{E}(\rho(z_{b+1}, M_t^{j+})|D).$$

Let $e_{j+1}$ be the edge that rings at time $\tau_{j+1}$ and let $\{\tilde{M}_t^j\}$ be the process that behaves exactly like $\{M_t^j\}$, except that if the balls at the endpoints of $e_{j+1}$ are red and white, it switches them if and only if $\{M_t^j\}$ *does not* switch them. Clearly $\tilde{M}_t^j$ has the same distribution as $M_t^j$ and hence

(8) $$\mathbf{E}(\rho(z_{b+1}, M_t^j)|D) = \mathbf{E}(\rho(z_{b+1}, \tilde{M}_t^j)|D).$$

But note that

$$\tfrac{1}{2}\rho(z_{b+1}, M_t^j) + \tfrac{1}{2}\rho(z_{b+1}, \tilde{M}_t^j) = \rho(z_{b+1}, M_t^{j+}).$$

Taking a conditional expectation given $D$ and combining with (8) gives (7). It remains to show that

(9) $$\mathbf{E}(\rho(z_{b+1}, M_t^{j+})|D) = \mathbf{E}(\rho(z_{b+1}, M_t^{j+1})|D).$$

Recall that $M^{j+1}$ behaves exactly like $M^{j+}$ except that it performs a depinking step at time $\tau_{j+1}$. Let $\widehat{M}^{j+1}$ be the process that behaves exactly like $M^{j+1}$, except that if it depinks at time $\tau_{j+1}$ then it makes the opposite recoloring choice (i.e., if $M^{j+1}$ recolors its pink balls white, then $\widehat{M}^{j+1}$ recolors its pink balls red and vice versa). Clearly $\widehat{M}^{j+1}$ has the same distribution as $M^{j+1}$ and hence

(10) $$\mathbf{E}(\rho(z_{b+1}, \widehat{M}_t^{j+1})|D) = \mathbf{E}(\rho(z_{b+1}, M_t^{j+1})|D).$$

But note that

$$\tfrac{1}{2}\rho(z_{b+1}, \widehat{M}_t^{j+1}) + \tfrac{1}{2}\rho(z_{b+1}, M_t^{j+1}) = \rho(z_{b+1}, M_t^{j+}).$$

Taking a conditional expectation given $D$ and combining with (10) gives (9). □

We are now ready for Lemma 2. Let $\mathcal{U}(\cdot|z_1,\ldots,z_b)$ denote the uniform distribution over $V - \{z_1,\ldots,z_b\}$.

LEMMA 2. *Fix $t > 0$ and let $(Z_1,\ldots,Z_b) = v_t(1,\ldots,b)$. Consider a chameleon process with $b$ black balls. Then*

(11) $$\mathbf{E}(\|\mu_t(\cdot|Z_1,\ldots,Z_b) - \mathcal{U}(\cdot|Z_1,\ldots,Z_b)\|) \leq \frac{1}{m}\widehat{\mathbf{E}}\left(w_t + \frac{1}{2}p_t\right),$$

*where $\|\cdot\|$ denotes variation distance and we write $\widehat{\mathbf{E}}(\cdot)$ for $\mathbf{E}(\cdot|A)$.*

Note that since $Z_1, \ldots, Z_b$ are random variables (whose joint distribution is governed by $\mu_t$), so is the quantity $\|\mu_t(\cdot|Z_1, \ldots, Z_b) - \mathcal{U}(\cdot|Z_1, \ldots, Z_b)\|$; the left-hand side of (11) is the expectation of this random variable.

PROOF OF LEMMA 2. Recall that for probability distributions $\nu_1, \nu_2$, the total variation distance

$$\|\nu_1 - \nu_2\| = \sup_Q(\nu_2(Q) - \nu_1(Q)).$$

Fix distinct vertices $z_1, \ldots, z_b$ and a subset $Q \subset V - \{z_1, \ldots, z_b\}$. Let

$$z = (z_1, \ldots, z_b); \qquad Z = (Z_1, \ldots, Z_b)$$

and $m = n - b$. Lemma 1 implies that

(12) $\quad \mathcal{U}(Q|z_1, \ldots, z_b) - \mu_t(Q|z_1, \ldots, z_b) = \dfrac{|Q|}{m} - \mathbf{E}(\rho_t(Q)|Z = z),$

where we write $\rho_t(Q)$ for $\sum_{x \in Q} \rho_t(x)$. But

$$\mathbf{E}(\rho_t(Q)|Z = z) \geq \frac{1}{m}\widehat{\mathbf{E}}(\rho_t(Q)|Z = z),$$

since $\mathbf{P}(A) = \frac{1}{m}$. Hence the right-hand side of (12) is at most

$$\frac{1}{m}\widehat{\mathbf{E}}(|Q| - \rho_t(Q)|Z = z) = \frac{1}{m}\widehat{\mathbf{E}}\bigg(|Q \cap W_t| + \frac{1}{2}|Q \cap P_t|\bigg|Z = z\bigg)$$
$$\leq \frac{1}{m}\widehat{\mathbf{E}}\bigg(w_t + \frac{1}{2}p_t\bigg|Z = z\bigg).$$

Since this holds for all $Q$, we have

$$\|\mu_t(\cdot|z_1, \ldots, z_b) - \mathcal{U}(\cdot|z_1, \ldots, z_b)\| \leq \frac{1}{m}\widehat{\mathbf{E}}\bigg(w_t + \frac{1}{2}p_t\bigg|Z = z\bigg)$$

and hence

(13) $\quad \mathbf{E}(\|\mu_t(\cdot|Z_1, \ldots, Z_b) - \mathcal{U}(\cdot|Z_1, \ldots, Z_b)\|) \leq \dfrac{1}{m}\mathbf{E}\bigg(\widehat{\mathbf{E}}\bigg(w_t + \dfrac{1}{2}p_t\bigg|Z\bigg)\bigg).$

Recall that the event $A$ is independent of the trajectories of the balls (and hence independent of $Z$). It follows that we can rewrite the RHS of (13) as $\frac{1}{m}\widehat{\mathbf{E}}(w_t + \frac{1}{2}p_t)$. □

For a chameleon configuration $M$, let $S(M) = (r(M) + \frac{1}{2}p(M))$ and let $s(M) = \frac{S(M)}{m}$. For a configuration $M$, let $\mathcal{W}(M) = \mathbf{E}_M(T_1)$, where $T_1$ is the first depinking time. So $\mathcal{W}_t$ is the expected waiting time until the next depinking after time $t$, given $M_t$. In Section 4, we show that

$$\mathcal{W}_t \leq g(S_t),$$

where $g(S) = c \log d(S+b)^{2/d}$ for a universal constant $c$. Let

$$\gamma = \left(\frac{2}{\sqrt{4/3} + \sqrt{2/3}} - 1\right).$$

Note that $0 < \gamma < \sqrt{2} - 1$. Let

$$Y(M) = 1 + \frac{\gamma \mathcal{W}(M)}{g(S(M))}.$$

Then $1 \leq Y \leq 1 + \gamma$. For $0 \leq u \leq 1$ define $u^\sharp = \min(u, 1-u)$. We have

$$\widehat{\mathbf{E}}(1 - s_t) \leq \widehat{\mathbf{E}}(1 - s_t)Y_t$$

$$\leq \widehat{\mathbf{E}}\left(\frac{\sqrt{s_t^\sharp}}{s_t} Y_t\right),$$

where the second inequality can be seen by considering the cases $s_t \geq \frac{1}{2}$ and $s_t < \frac{1}{2}$ separately. Define $Z(M) = \frac{\sqrt{s(M)^\sharp}}{s(M)} Y(M)$ and let $L_t = \widehat{\mathbf{E}}(Z_t)$, so that $\widehat{\mathbf{E}}(1 - s_t) \leq L_t$. Note that

$$L_t = \sum_M \mathbf{P}(M_t = M) Z(M)$$

is differentiable. Define $f$ by

$$f(z) = \begin{cases} \frac{1}{4} g((1+\gamma)^2 m/z^2)^{-1}, & \text{if } z \geq \sqrt{2}, \\ \frac{1}{4}(c \log d)^{-1} L^{-2}, & \text{otherwise.} \end{cases}$$

Note that $g$ is increasing and, for $z \geq \sqrt{2}$, we have $g((1+\gamma)^2 m/z^2) \leq g(m) = c \log dL^2$ [recall that $(\gamma+1)^2 \leq 2$]. It follows that $f$ is increasing. For $Z_t < \sqrt{2}$, we have

$$f(Z_t) = \tfrac{1}{4}(c \log d)^{-1} L^{-2} \leq \tfrac{1}{4}(c \log d(S_t + b))^{-2/d} = \tfrac{1}{4} g(S_t)^{-1},$$

and, when $Z_t \geq \sqrt{2}$, we have $s_t = s_t^\sharp$ and hence

$$f(Z_t) = \frac{1}{4} g\left(\frac{(1+\gamma)^2}{Y_t^2} S_t\right)^{-1} \leq \frac{1}{4} g(S_t)^{-1}.$$

Thus $f(Z_t) \leq \frac{1}{4g(S_t)}$ for all $t$.

LEMMA 3. *The differential equation*

(14) $$L'_t \leq -\gamma L_t f(L_t/2)$$

*is satisfied.*

PROOF. First, we will show that for every $t \geq 0$,

$$\widehat{\mathbf{E}}(Z_{t+\varepsilon}|M_t) \leq Z_t(1 - 2\gamma\varepsilon f(Z_t)) + O(\varepsilon^2). \tag{15}$$

By the Markov property, it is enough to show that for any starting configuration $M_0$ we have

$$\widehat{\mathbf{E}}_{M_0}(Z_\varepsilon) \leq Z_0(1 - 2\gamma\varepsilon f(Z_0)) + O(\varepsilon^2). \tag{16}$$

Note that $s_t$ is constant on each interval $[T_j, T_{j+1})$ and that $s_{T_j} = \frac{1}{m}r_{T_j}$ for all $j$. It follows that $\{s_t : t \geq 0\}$ is a martingale and hence $\mathbf{P}(A|s_0 = s) = \frac{s}{m}$ for every $s$. Using this, it is easy to check that

$$\widehat{\mathbf{E}}(h(M_t)) = \mathbf{E}\left(\frac{s_t}{s_0}h(M_t)\right) \tag{17}$$

for every function $h$ and $t \geq 0$. Thus

$$\widehat{\mathbf{E}}(Z_\varepsilon) = \mathbf{E}\left(Z_\varepsilon \frac{s_\varepsilon}{s_0}\right) = \frac{1}{s_0}\mathbf{E}(\sqrt{s_\varepsilon^\sharp}Y_\varepsilon). \tag{18}$$

We claim that

$$\mathbf{E}(\sqrt{s_\varepsilon^\sharp}Y_\varepsilon) \leq \left(1 - \frac{\gamma\varepsilon}{2g(S_0)}\right)\sqrt{s_0^\sharp}Y_0 + O(\varepsilon^2). \tag{19}$$

Combining this with (18) gives

$$\widehat{\mathbf{E}}(Z_\varepsilon) \leq \left(1 - \frac{\gamma\varepsilon}{2g(S_0)}\right)\frac{\sqrt{s_0^\sharp}}{s_0}Y_0 + O(\varepsilon^2)$$
$$\leq (1 - 2\gamma\varepsilon f(Z_0))Z_0 + O(\varepsilon^2),$$

which verifies (16).

It remains to verify (19). Let $\{\tilde{M}_t : t \geq 0\}$ be the process that behaves like $M_t$, except that it does not perform any depinking steps. So $\tilde{M}_t$ pinkens, but does not depink and hence accumulates pink balls. Let $\tilde{p}(\cdot), \tilde{r}(\cdot)$ and $\tilde{w}(\cdot)$ output the number of pink, red and white balls, respectively, in a configuration $\tilde{M}$, and define

$$\mathcal{W}(\tilde{M}) = \mathbf{E}_{\tilde{M}}(\inf\{t : \tilde{p}_t \geq \min(\tilde{r}_t, \tilde{w}_t)\}),$$

where we write $\mathbf{E}_{\tilde{M}}(\cdot) := \mathbf{E}(\cdot|\tilde{M}_0 = \tilde{M})$. Let $T_1$ be the first time that $M_t$ depinks. If we start with $\tilde{M}_0 = M_0$, then for every $t \geq 0$ we have

$$\mathcal{W}(\tilde{M}_t) = \begin{cases} \mathcal{W}(M_t), & \text{if } t < T_1; \\ 0, & \text{if } t \geq T_1. \end{cases}$$

By analogy with the previous discussion, let $\tilde{S} = (\tilde{w} + \frac{1}{2}\tilde{p})$, let $\tilde{s} = \frac{1}{m}\tilde{S}$ and define $\tilde{Y}(\tilde{M}) = 1 + \frac{\gamma \mathcal{W}(\tilde{M})}{g(\tilde{S})}$. Fix a starting configuration $\tilde{M}_0 = M_0$ with $\tilde{p} \leq \min(\tilde{r}, \tilde{w})$. We will show that

$$\mathbf{E}(\sqrt{\tilde{s}_\varepsilon^\sharp} \tilde{Y}_\varepsilon) \leq \left(1 - \frac{\gamma \varepsilon}{2g(\tilde{S}_0)}\right)\sqrt{\tilde{s}_0^\sharp}\tilde{Y}_0 + O(\varepsilon^2) \tag{20}$$

and

$$\mathbf{E}(\sqrt{s_\varepsilon}Y_\varepsilon) \leq \mathbf{E}(\sqrt{\tilde{s}_\varepsilon}\tilde{Y}_\varepsilon) + O(\varepsilon^2). \tag{21}$$

These will imply (19), since $\tilde{Y}_0 = Y_0$ and $\tilde{s}_0 = s_0$. Note that since $\tilde{M}_t$ never depinks, $\tilde{s}_t \equiv \tilde{s}$ is constant and we have

$$\mathbf{E}(\sqrt{\tilde{s}_\varepsilon^\sharp}\tilde{Y}_\varepsilon) - \mathbf{E}(\sqrt{\tilde{s}_0^\sharp}\tilde{Y}_0) = \sqrt{\tilde{s}^\sharp} \frac{\gamma \mathbf{E}(\mathcal{W}(\tilde{M}_\varepsilon) - \mathcal{W}(\tilde{M}_0))}{g(\tilde{S})}. \tag{22}$$

Let $B = \{\tilde{M} : \tilde{p} \geq \min(\tilde{r}, \tilde{w})\}$. For any configuration $\tilde{M}$, the quantity $\mathcal{W}(\tilde{M})$ is the expected waiting time to hit $B$ starting from $\tilde{M}$. Hence Proposition 11 below implies that $\mathbf{E}(\mathcal{W}(\tilde{M}_\varepsilon) - \mathcal{W}(\tilde{M}_0)) = -\varepsilon + O(\varepsilon^2)$. Combining this with the fact that $\tilde{Y}_0 \leq 2$ gives

$$\mathbf{E}(\sqrt{\tilde{s}_\varepsilon^\sharp}Y_\varepsilon) - \mathbf{E}(\sqrt{\tilde{s}_0^\sharp}Y_0) \leq \sqrt{\tilde{s}^\sharp}\frac{-\gamma \varepsilon}{2g(\tilde{S})}\tilde{Y}_0 + O(\varepsilon^2) \tag{23}$$

and (20) is verified. It remains to verify (21). Let $D$ be the event that there is exactly one transition in $[0, \varepsilon]$ and let $F$ be the event that $\tilde{M}_\varepsilon \in B$. Note that $\mathbf{P}(F^c \cup D) \geq 1 - O(\varepsilon^2)$, and on the event $F^c$ we have $\tilde{M}_\varepsilon = M_\varepsilon$ and $\tilde{Y}_\varepsilon = Y_\varepsilon$. Thus, it is enough to show that

$$\mathbf{E}(\sqrt{s_\varepsilon^\sharp}Y_\varepsilon | D \cap F) \leq \mathbf{E}(\sqrt{\tilde{s}_\varepsilon^\sharp}\tilde{Y}_\varepsilon | D \cap F). \tag{24}$$

On the event $D \cap F$, $M_\varepsilon$ is obtained from $\tilde{M}_\varepsilon$ by performing a depinking operation. Fix $\tilde{M} \in B$ with $\tilde{s} \leq \frac{1}{2}$ and let $M$ be the configuration obtained from $\tilde{M}$ by depinking. Then $M$ is a random configuration with two possible values, since the pink balls can be recolored either red or white. Let $Y = Y(M)$ and $s = s(M)$. Since $Y \leq 1 + \gamma$ and $s^\sharp \leq s$, we have

$$\mathbf{E}(\sqrt{s^\sharp}Y) \leq (1 + \gamma)\mathbf{E}(\sqrt{s})$$
$$= (1 + \gamma)\left(\frac{1}{2}\sqrt{\tilde{s} + \Delta(\tilde{s})} + \frac{1}{2}\sqrt{\tilde{s} - \Delta(\tilde{s})}\right)$$
$$\leq (1 + \gamma)\left(\frac{\sqrt{4/3} + \sqrt{2/3}}{2}\right)\sqrt{\tilde{s}} = \sqrt{\tilde{s}},$$

where the second inequality holds because $\Delta(\tilde{s}) \geq \tilde{s}/3$. But, $\tilde{s} = \tilde{s}^\sharp$ and $\tilde{Y} \geq 1$, so $\mathbf{E}(\sqrt{s^\sharp}Y) \leq \sqrt{\tilde{s}^\sharp}\tilde{Y}$. This holds for all $\tilde{M} \in B$ with $\tilde{s} \leq 1/2$. By reversing the roles of red and white balls, we can argue similarly in the case $\tilde{s} > \frac{1}{2}$. This verifies (24), which in turn verifies (15).

Taking the expectation (conditional on $A$) of both sides of (15), we get

$$L_{t+\varepsilon} \leq L_t - 2\gamma\varepsilon\widehat{\mathbf{E}}(Z_t f(Z_t)) + O(\varepsilon^2). \tag{25}$$

We will use the following lemma from [10].

LEMMA 4 ([10]). *Suppose that $Z \geq 0$ is a nonnegative random variable and $f$ is a nonnegative increasing function. Then $\mathbf{E}(Zf(2Z)) \geq \frac{\mathbf{E}Z}{2} \cdot f(\mathbf{E}Z)$.*

Lemma 4 implies that $\widehat{\mathbf{E}}(Z_t f(Z_t)) \geq \frac{L_t}{2} f(\frac{L_t}{2})$. Plugging this into (25), we get

$$L_{t+\varepsilon} \leq L_t - \varepsilon\gamma L_t f(L_t/2) + O(\varepsilon^2) \tag{26}$$

and the proof is complete. □

REMARK. Note that since $f$ is bounded below by $\Omega((\log dL^2)^{-1})$, Lemmas 2, 3 and 12 imply that the spectral gap for the simple exclusion process must also be at least of order $(\log dL^2)^{-1}$.

We are now ready to prove our final ingredient, namely, a bound on the RHS of (11) when $G$ is the torus:

LEMMA 5. *Let $G = \mathbf{Z}^d/L\mathbf{Z}^d$. Fix $b \in \{0, \ldots, k-1\}$ and consider a chameleon process on $G$ with $b$ black balls. There is universal constant $C$ such that if $t \geq C[d \log dL^2 + \log dL^2 \log k]$, then*

$$\frac{1}{m}\widehat{\mathbf{E}}\left(w_t + \frac{1}{2}p_t\right) \leq \frac{1}{4k}. \tag{27}$$

PROOF. We will need the following lemma from [10].

LEMMA 6 ([10]). *Suppose that $L_t$ is a nonnegative, differentiable function of $t$ that satisfies the differential equation $L'_t \leq -\gamma L_t f(L_t/2)$. Then for every*

$$t \geq \int_{\delta/2}^{L_0/2} \frac{dz}{\gamma z f(z)},$$

*we have $L_t \leq \delta$.*

Recall that the number of vertices $n$ is $L^d$. Since $L_0 = Z_0 = \frac{1}{\sqrt{s_0}}Y_0 \leq 2\sqrt{n}$, Lemmas 6 and 3 imply that $L_t < \frac{1}{4k}$ for

$$t \geq \int_{1/8k}^{\sqrt{n}} \frac{dz}{\gamma z f(z)}. \tag{28}$$

We need to show that this integral is at most

(29) $$Cd\log d L^2 + C\log d L^2 \log k$$

for a constant $C$. We can write the integral as

(30) $$\int_{1/8k}^{\sqrt{2}} \frac{dz}{\gamma z f(z)} + \int_{\sqrt{2}}^{\sqrt{m/b}} \frac{dz}{\gamma z f(z)} + \int_{\sqrt{m/b}}^{\sqrt{n}} \frac{dz}{\gamma z f(z)} \equiv I_1 + I_2 + I_3.$$

We will bound the terms $I_1, I_2$ and $I_3$ separately. Since $\frac{1}{\gamma f(z)} = 4c\gamma^{-1}\log d L^2$ when $z < \sqrt{2}$,

$$I_1 = 4c\gamma^{-1}\log d L^2 \log(8\sqrt{2}k) \leq \alpha \log d L^2 \log k$$

for a constant $\alpha$ and hence the second term in (29) dominates $I_1$ if $C$ is sufficiently large. To bound $I_2$ and $I_3$, note that for $z \geq \sqrt{2}$ we have

(31) $$\frac{1}{\gamma f(z)} = 4c\gamma^{-1}\log d\left(\frac{m(1+\gamma)^2}{z^2} + b\right)^{2/d} \leq \beta \log d \cdot \max\left(\frac{m}{z^2}, b\right)^{2/d}$$

for a constant $\beta$. Hence

$$I_2 \leq \int_{\sqrt{2}}^{\sqrt{m/b}} \beta \log d \cdot \left(\frac{m}{z^2}\right)^{2/d} \frac{1}{z} dz$$
$$\leq \beta \log d\, m^{2/d}\left(\int_{\sqrt{2}}^{\infty} z^{-(1+4/d)}\, dz\right).$$

But, $m^{2/d} \leq n^{2/d} = L^2$ and an elementary calculation shows that the integral in the second line is at most $d$. Hence the first term in (29) dominates $I_2$ if $C$ is sufficiently large. It remains to bound $I_3$. Using (31) we get

(32) $$I_3 \leq \int_{\sqrt{m/b}}^{\sqrt{n}} \beta \log d \cdot b^{2/d} \frac{1}{z} dz$$
$$\leq \beta \log d\, b^{2/d}\left[\frac{1}{2}\log(n/m) + \frac{1}{2}\log b\right].$$

But since

$$m \geq n/2; \qquad b^{2/d} \leq n^{2/d} = L^2; \qquad b \leq k,$$

we have $I_3 = O(L^2 \log d \log k)$. Thus $I_3$ is dominated by the second term in (29) if $C$ is sufficiently large.

We have shown that $I_1$, $I_2$ and $I_3$ are all of the correct form. This completes the proof of Lemma 5. $\square$

We are now ready for our main theorem.

THEOREM 7. *The mixing time $\tau_{\mathrm{mix}}$ for $k$-particle exclusion in $\mathbf{Z}^d/L\mathbf{Z}^d$ satisfies*

$$\tau_{\mathrm{mix}} \leq C[d \log dL^2 + \log dL^2 \log k] \tag{33}$$

*for a universal constant $C$.*

PROOF. Let $G = \mathbf{Z}^d/L\mathbf{Z}^d$. Fix $t \geq 0$ and suppose that $(Z_1, \ldots, Z_k) \sim \mu_t$. Lemma 12 in Section 5.3 gives

$$\|\mu_t - \mathcal{U}\| \leq \sum_{b=0}^{k-1} \mathbf{E}(\|\mu_t(\cdot|Z_1, \ldots, Z_b) - \mathcal{U}(\cdot|Z_1, \ldots, Z_b)\|), \tag{34}$$

but by Lemmas 2 and 5, each term in the RHS of (34) is at most $\frac{1}{4k}$ whenever $t \geq C[d \log dL^2 + \log dL^2 \log k]$ for a universal constant $C$. Hence the sum is at most $\frac{1}{4}$. □

**4. Waiting to depink.** Let $G = \mathbf{Z}^d/L\mathbf{Z}^d$. Consider a chameleon process on $G$ with $b$ black balls. Let $T_1$ be the time of the first depinking. The principal result of this section is an upper bound on $\mathbf{E}(T_1)$, valid uniformly over starting states $M_0$.

Recall that $W_t$ and $R_t$ denote the set of balls that are white and red, respectively, at time $t$. Say an edge is *conflicting at time $t$* if one of its endpoints contains a ball from $W_0$ and the other contains a ball from $R_0$.

We will need the following lemma.

LEMMA 8. *There exist universal positive constants $c$ and $B$, with $B \leq 2/3$, such that if $\tau \sim \mathrm{uniform}[0, c \log d(b+S)^{2/d}]$, then the expected number of edges conflicting at time $\tau$ is at least $Br_0 d$.*

PROOF. Let $G = (V, E)$ be a connected graph. Define the *random walk on $G$ with rate $\lambda$* as the continuous-time Markov chain on $V$ with the following transition rule: If the current state is $x \in V$, traverse each edge $e$ incident to $x$ at rate $\lambda$. Let $G = \mathbf{Z}^d/L\mathbf{Z}^d$ and let $\{Z_t : t \geq 0\}$ be the random walk on $G$ with rate 2. Let $G'$ be the subgraph of $G$ induced by removing the origin and let $\{Z'_t : t \geq 0\}$ be the random walk on $G'$ with rate 2.

Consider the chameleon (or interchange) process on $G$. Clearly, for any distinct pair of balls $x$ and $y$, we have

$$\mathbf{P}(v_t(x) \text{ is adjacent to } v_t(y)) = \mathbf{P}(Z'_t \in \mathcal{B}|Z'_0 = v_0(x) - v_0(y)), \tag{35}$$

where $\mathcal{B}$ is the set of neighbors of the origin in $G$.

The Markov chain $\{Z_t\}$ is easier to work with than $\{Z'_t\}$. They are related in the following way: Fix $z \in V - \{0\}$ and $\alpha > 0$, and let $\tau \sim \mathrm{uniform}[0, \alpha]$. Then

$$\mathbf{P}(Z'_\tau \in \mathcal{B}|Z'_0 = z) \geq \mathbf{P}(Z_\tau \in \mathcal{B}|Z_0 = z). \tag{36}$$

The reason that (36) holds is that the "wasted time" $Z_t$ spends in 0 lowers the probability that $Z_\tau \in \mathcal{B}$. To give a formal justification for (36), we will show that $|\{t < \alpha : Z'_t \in \mathcal{B}\}|$ stochastically dominates $|\{t < \alpha : Z_t \in \mathcal{B}\}|$. Let $\tilde{Z}_t$ be the process obtained from $Z_t$ by "skipping past" all visits to the origin, that is, the process obtained by replacing transitions to 0 by transitions to a uniformly chosen element of $\mathcal{B}$.

Note that in the natural coupling (where the trajectory of $\{\tilde{Z}_t\}$ is constructed from that of $\{Z_t\}$ by simply deleting every visit to the origin), we have

$$|\{t < \alpha : \tilde{Z}_t \in \mathcal{B}\}| \geq |\{t < \alpha : Z_t \in \mathcal{B}\}|.$$

Also, clearly we can couple $(\tilde{Z}_t, Z'_t)$ so that

$$\{t : \tilde{Z}_t \in \mathcal{B}\} = \{t : Z'_t \in \mathcal{B}\}.$$

It follows that $|\{t < \alpha : Z'_t \in \mathcal{B}\}|$ stochastically dominates $|\{t < \alpha : Z_t \in \mathcal{B}\}|$ and (36) is verified.

Let $p(\cdot, \cdot)$ be the transition kernel for $\{Z_t\}$ and, for $\varepsilon \geq \frac{1}{n}$, define

$$\tau(\varepsilon) = \inf\{t : p^t(x,y) \leq \tfrac{7}{6}\varepsilon \ \forall x,y\}.$$

We will need the following easy proposition, which is proved in the next section.

PROPOSITION 9. *For any $\varepsilon \geq \frac{1}{n}$, we have*

$$\tau(\varepsilon) \leq \varepsilon^{-2/d} D \log d,$$

*where $D$ is a universal constant.*

We are now ready to complete our proof of Lemma 8. Let $c = 20D$, where $D$ is the constant that appears in Proposition 9, and let $B = 1/72$. For $x \in R_0$, let $N_x = \{y \in W_0 : v_\tau(y) \text{ is adjacent to } v_\tau(x)\}$. It is enough to show that for every $x \in R_0$, we have $\mathbf{E}(N_x) \geq Bd$. But

$$\begin{align}
(37) \quad \mathbf{E}(N_x) &= \sum_{y \in W_0} \mathbf{P}(v_\tau(x) \text{ is adjacent to } v_\tau(y)) \\
(38) \quad &= \sum_{y \in W_0} \mathbf{P}(Z'_\tau \in \mathcal{B} | Z'_0 = v_0(x) - v_0(y)) \\
(39) \quad &\geq \sum_{y \in W_0} \mathbf{P}(Z_\tau \in \mathcal{B} | Z_0 = v_0(x) - v_0(y)) \\
(40) \quad &= \mathbf{E}(p^\tau(\mathcal{A}, \mathcal{B})),
\end{align}$$

where $\mathcal{A} = \{v_0(x) - v_0(y) : y \in W_0\}$. The second line follows from (35) and the inequality follows from (36), but the symmetry of $\{Z_t\}$ implies that the quantity (40) is

$$\mathbf{E}(p^\tau(\mathcal{B}, \mathcal{A})) = |\mathcal{B}|\mathbf{P}_{\mathcal{U}_B}(Z_\tau \in \mathcal{A})$$
$$\geq d \min_{x \in V} \mathbf{P}_x(Z_\tau \in \mathcal{A}),$$

where $\mathcal{U}_B$ denotes the uniform distribution over $\mathcal{B}$.

We will now complete the proof by showing that $\mathbf{P}_x(Z_\tau \in \mathcal{A}) \geq 1/72$ for all $x \in V$. Fix $x \in V$. Since $\tau \sim \text{uniform}[0, c \log d(b+S)^{2/d}]$ and $c = 20D$, we have

$$\mathbf{P}_x(Z_\tau \in \mathcal{A}) \geq \tfrac{1}{2}\mathbf{P}_x(Z_\tau \in \mathcal{A} | \tau \geq 10D \log d(b+S)^{2/d}).$$

Thus, it is enough to show that $p^t(x, \mathcal{A}) \geq 1/36$ for $t > 10D \log d(b+S)^{2/d}$. Let $\mathcal{A}^c$ denote the complement of $\mathcal{A}$ in $V$. Then

$$p^t(x, \mathcal{A}) = 1 - p^t(x, \mathcal{A}^c)$$
(41)
$$\geq 1 - |\mathcal{A}^c| \max_{y \in V} p^t(x, y)$$
$$= 1 - (b + r_0 + p_0) \max_{y \in V} p^t(x, y).$$

Since $b \leq n/2$ and $w_0 \geq r_0 \geq p_0$, we have $w_0 \geq n/6$ and hence $(b + r_0 + p_0) \leq 5n/6$. Plugging $\varepsilon = \tfrac{5}{6}\tfrac{1}{(b+r_0+p_0)}$ into Proposition 9 shows that for $t \geq D \log d \cdot [\tfrac{6}{5}(b + r_0 + p_0)]^{2/d}$, we have

$$p^t(x, y) \leq \frac{35}{36} \frac{1}{(b + r_0 + p_0)}$$

for any $y$ and hence $p^t(x, \mathcal{A}) \geq 1/36$ by (41). Finally, note that $S \geq (r_0 + p_0)/2$ and hence

$$D \log d \cdot [\tfrac{6}{5}(b + r_0 + p_0)]^{2/d} \leq D \log d[3(b+S)]^{2/d}$$
$$\leq 10D \log d(b+S)^{2/d}. \qquad \square$$

We are now ready for the main result of this section. Recall that $T_1$ denotes the time of the first depinking.

LEMMA 10. *Let $S = r_0 + \tfrac{1}{2}p_0$ and suppose $S \notin \{0, m\}$. Then*

(42) $$\mathbf{E}(T_1) \leq c \log d(b+S)^{2/d}$$

*for a universal constant $c$.*

PROOF. Since the chameleon process is symmetric with respect to the roles of the red and white balls, and since the expression on the RHS of (42) is increasing in $S$, we may assume that $S \leq m/2$. Note that $r_t + \frac{1}{2}p_t = S$ for all $t \in [0, T_1)$. For $\beta > 0$, let

$$T_\beta = \min\{t : \text{at least } \beta r_0 \text{ balls have been pinkened by time } t\}.$$

We will show that there exist universal constants $c, \beta > 0$ such that

(43) $$\mathbf{P}(T_\beta \leq c \log d(b+S)^{2/d}) \geq \beta.$$

Note that the lemma follows from (43) since we can incorporate an extra factor $\lceil \frac{1}{2\beta^2} \rceil$ into $c$.

It remains to verify (43). Recall that $W_t$ and $R_t$ denote the set of balls that are white and red, respectively, at time $t$. Say edges $e$ and $f$ are neighbors if they are incident to a common vertex. Note that for any edge $e$,

(44) $$\mathbf{P}(e \text{ rings for some time } t \in [0,1] \text{ and before any of its neighbors}) \geq \frac{1}{5d}.$$

To see this, fix an edge $e$. Let $A_1$ be the event that $e$ rings before any of its neighbors and let $A_2$ be the event that $e$ rings for some $t \in [0,1]$. Since no edge in $\mathbf{Z}^d/L\mathbf{Z}^d$ has more than $4d-2$ neighbors, we have $\mathbf{P}(A_1) \geq \frac{1}{4d-1}$, and since the edges ring at rate 2, we have $\mathbf{P}(A_2) = 1 - e^{-2}$. Finally, note that $A_1$ and $A_2$ are positively correlated and hence

$$\mathbf{P}(A_1 A_2) \geq \mathbf{P}(A_1)\mathbf{P}(A_2) \geq \frac{1-e^{-2}}{4d-1} \geq \frac{1}{5d}.$$

Recall that an edge is conflicting at time $t$ if one of its endpoints contains a ball from $W_0$ and the other contains a ball from $R_0$. Let $N$ denote the number of edges in $G$ that:

1. Conflict at time $\tau$.
2. Ring before any of their neighbors after time $\tau$.
3. Ring for some $t \in [\tau, \tau+1]$.

Lemma 8 says that the expected number of edges that conflict at time $\tau$ is at least $Br_0 d$ for a constant $B$. Combining this with (44) and using the strong Markov property, we get

$$\mathbf{E}(N) \geq [Br_0 d]\frac{1}{5d} = \frac{Br_0}{5}.$$

Let $J$ be the number of balls pinkened before time $T_1 \wedge (\tau+1)$. Note that $J \geq N$. Thus $\mathbf{E}(J) \geq \mathbf{E}(N) \geq Br_0/5$. Hence, an application of Markov's inequality to the nonnegative random variable $2r_0 - J$ gives

$$\mathbf{P}\left(J \geq \frac{Br_0}{10}\right) \geq \frac{B}{20-B}$$
$$> \frac{B}{20}.$$

Finally, note that on the event that $J \geq Br_0/10$, we have $T_{B/10} \leq \tau + 1$. It follows that

$$\mathbf{P}(T_{B/10} \leq c \log d(S+b)^{2/d} + 1) > B/20.$$

Incorporating the (small) additive constant $1/\log d(S+b)^{2/d}$ into $c$ and setting $\beta = B/20$ verifies (43). □

**5. Technical details.** We now prove three technical results that were needed in previous sections.

5.1. *Easy proposition used in the proof of Lemma 3.*

PROPOSITION 11. *Let $X_t$ be a continuous-time Markov chain on a finite state space $V$ and let $B \subset V$ be a subset of states. For $x \in V$, let $\mathcal{W}(x) = \mathbf{E}_x(\inf\{t : X_t \in B\})$ be the expected time to hit $B$ starting from $x$. Then, for every $x \in V - B$, we have*

$$\mathbf{E}_x(\mathcal{W}(X_\varepsilon) - \mathcal{W}(x)) = -\varepsilon + O(\varepsilon^2).$$

PROOF. Fix $x \in V - B$ and let $X_0 = x$. Let $q_x$ be the rate of transitions out of $x$, let $Y$ be the next state visited after $x$ and let $W = \inf\{t : X_t = Y\}$. Then $W$ has an exponential distribution with parameter $q_x$. It follows that $\mathcal{W}(x) = 1/q_x + \mathbf{E}(\mathcal{W}(Y))$. Let $D$ be the event that the chain makes exactly one transition in $[0, \varepsilon]$. Then

$$\mathbf{E}(\mathcal{W}(X_\varepsilon) - \mathcal{W}(x)|D) = \mathbf{E}(\mathcal{W}(Y)) - \mathcal{W}(x) = -1/q_x$$

and $\mathbf{P}(D) = \varepsilon q_x + O(\varepsilon^2)$. Hence

$$\mathbf{E}(\mathcal{W}(X_\varepsilon) - \mathcal{W}(x)) = (\varepsilon q_x) \cdot (-1/q_x) + O(\varepsilon^2)$$
$$= -\varepsilon + O(\varepsilon^2). \qquad \square$$

5.2. *Probability bounds for random walk on the torus.*

PROOF OF PROPOSITION 9. We adopt the notation of the previous section. Let $X_t$ be the random walk of rate 2 on the cycle $C_L$ of length $L$. It is well known that for all $j \leq L, l \geq 0$ and $t \geq j^2 l$, we have

$$\mathbf{P}(X_t = y) \leq \frac{1}{j}(1 + e^{-\gamma l})$$

for every $y \in C_L$, where $\gamma > 0$ is a universal constant. Thus the proposition holds when $d = 1$ (provided that we replace $\log d$ by 1 in the bound).

When $d > 1$, the coordinates of $Z_t$ are independent copies of $X_t$. Thus, for every $z \in V$ and $t \geq j^2 l$, we have

$$\mathbf{P}(Z_t = z) \leq \frac{1}{j^d}(1 + e^{-\gamma l})^d.$$

Plugging in $l = \frac{4}{\gamma} \log d$ we get

$$\mathbf{P}(Z_t = z) \leq \frac{1}{j^d}\left(1 + \frac{1}{d^4}\right)^d \leq \frac{1}{j^d} \times \frac{7}{6},$$

where we have used the inequality $(1 + 1/d^4)^d \leq \frac{7}{6}$, valid for $d \geq 2$. The above shows that for all $j \leq n$, we have

$$\tau\left(\frac{1}{j^d}\right) \leq j^2 \frac{4}{\gamma} \log d.$$

Setting $\varepsilon = 1/j^d$ and $D = 4/\gamma$ we get

$$\tau(\varepsilon) \leq \varepsilon^{-2/d} D \log d. \qquad \square$$

5.3. *Lemma about total variation distance.*

LEMMA 12. *Fix a graph $G$, let $\mu$ and $\nu$ be arbitrary probability distributions supported on $k$-tuples of distinct vertices of $G$ and suppose that $(Z_1, \ldots, Z_k) \sim \mu$. Then*

$$(45) \quad \|\mu - \nu\| \leq \sum_{l=0}^{k-1} \mathbf{E}(\|\mu(\cdot|Z_1, \ldots, Z_l) - \nu(\cdot|Z_1, \ldots, Z_l)\|).$$

PROOF. Recall that for probability distributions $\widehat{\mu}$ and $\widehat{\nu}$, the total variation distance

$$(46) \quad \|\widehat{\mu} - \widehat{\nu}\| = \min_{W_1 \sim \widehat{\mu}, W_2 \sim \widehat{\nu}} \mathbf{P}(W_1 \neq W_2).$$

Thus, for every $l$ and $z_1, \ldots, z_l$, we can construct $W_1 \sim \mu(\cdot|z_1, \ldots, z_l)$ and $W_2 \sim \nu(\cdot|z_1, \ldots, z_l)$ such that

$$\mathbf{P}(W_1 \neq W_2) = \|\mu(\cdot|z_1, \ldots, z_l) - \nu(\cdot|z_1, \ldots, z_l)\|.$$

We couple $Z \sim \mu$ with $Y \sim \nu$ as follows. Choose $(Y_1, Z_1)$ according to the optimal coupling [i.e., a coupling that achieves the minimum in the RHS of (46)] and subsequently, for all $l$ with $1 \leq l \leq k - 1$: if $(Z_1, \ldots, Z_l) = (Y_1, \ldots, Y_l)$, then choose $(Z_{l+1}, Y_{l+1})$ according to the optimal coupling of $\mu(\cdot|Z_1, \ldots, Z_l)$ and $\nu(\cdot|Z_1, \ldots, Z_l)$; else couple $(Y_{l+1}, Z_{l+1})$ in an arbitrary way.

Note that

$$\mathbf{P}(Z \neq Y) = \sum_{l=0}^{k-1} \mathbf{P}((Z_1, \ldots, Z_l) = (Y_1, \ldots, Y_l), Z_{l+1} \neq Y_{l+1}), \tag{47}$$

but on the event that $(Z_1, \ldots, Z_l) = (Y_1, \ldots, Y_l)$, the pair $(Z_{l+1}, Y_{l+1})$ is chosen according to the optimal coupling of $\mu(\cdot|Z_1, \ldots, Z_l)$ and $\nu(\cdot|Z_1, \ldots, Z_l)$. Hence the RHS of (47) is at most

$$\sum_{l=0}^{k-1} \mathbf{E}(\|\mu(\cdot|Z_1, \ldots, Z_l) - \nu(\cdot|Z_1, \ldots, Z_l)\|). \qquad \square$$

**6. Lower bounds.** In this section, we prove a matching lower bound of $\Omega(L^2 \log k)$ for the mixing time of the exclusion process in $\mathbf{Z}^d/L\mathbf{Z}^d$ with $k$ black balls. We will accomplish this by showing that there is a universal constant $c$ and an exchangeable event $A$, such that

$$\mathcal{U}(A) - \mu_t(A) > \tfrac{1}{4}$$

for every $t < cL^2 \log k$, where $\mathcal{U}$ and $\mu_t$ are defined as in Section 3.

The idea behind the proof is the following: If the black balls start out with their first coordinate close to $L/4$, then with high probability most will still be in the left half of the torus at time $t = \Omega(L^2 \log k)$. For convenience, we assume that $L$ is a multiple of 8. The proof for the general case is similar. Suppose that balls $1, \ldots, k$ are black and the remainder are white. Suppose that $v_t(x) = (v_t^1(x), \ldots, v_t^d(x))$, and let the starting configuration be chosen to minimize the quantity

$$\sum_{x=1}^{k} \left| v_0^1(x) - \frac{L}{4} \right|, \tag{48}$$

and such that $v_0^1(x) < L/2$ for every black ball $x$ (so that the first coordinate takes the value 0 in preference to $L/2$). For $v \in V$, suppose that $v = (v^1, \ldots, v^d)$ and let

$$S = \left\{ v \in V : \frac{L}{8} \leq v^1 \leq \frac{3L}{8} \right\}; \qquad U = \left\{ v \in V : \frac{L}{2} \leq v^1 < L \right\}.$$

Note that for every ball $x$, the process $\{v_t^1(x) : t \geq 0\}$ is just the continuous-time random walk of rate 2 on the cycle $C_L$ of length $L$. Hence $\mathbf{P}(v_t(x) \in U) \leq \tfrac{1}{2}$ for every black ball $x$ and $t \geq 0$. It is also well known that if $v_0(x) \in S$, then we have $\mathbf{P}(v_t(x) \in U) \leq \tfrac{1}{2} - \delta \exp(-\gamma t/L^2)$ for universal constants $\delta, \gamma > 0$. Let

$$N_t = |\{1 \leq x \leq k : v_t(x) \in U\}|$$

be the number of black balls in $U$ at time $t$. Since at least half of the black balls must start in $S$, we have $\mathbf{E}(N_t) \leq \frac{k}{2}[1 - \delta \exp(-\gamma t/L^2)]$ for every $t$.

For vertices $v \in U$, let $X^t(v)$ be the indicator for the event that $v$ is occupied by a black ball at time $t$, so that $N_t = \sum_{v \in U} X^t(v)$. Liggett [8] proved that for $u \neq v$, the random variables $X^t(u)$ and $X^t(v)$ are negatively correlated, that is, $\mathbf{E}(X^t(u)X^t(v)) \leq \mathbf{E}(X^t(u))\mathbf{E}(X^t(v))$. It follows that

$$\mathbf{E}(N_t^2) = \mathbf{E}\bigg(\sum_v X^t(v) + \sum_{u \neq v} X^t(u)X^t(v)\bigg)$$

$$\leq \mathbf{E}\bigg(\sum_v X^t(v)\bigg) + \sum_{u \neq v} \mathbf{E}(X^t(u))\mathbf{E}(X^t(v))$$

$$\leq \mathbf{E}(N_t) + (\mathbf{E}(N_t))^2$$

and hence $\text{var}(N_t) \leq \mathbf{E}(N_t) \leq k/2$. Thus Chebyshev's inequality gives

$$\mathbf{P}\bigg(N_t \geq \frac{k}{2}\bigg) \leq \frac{\text{var}(N_t)}{((k\delta/2)e^{-\gamma t/L^2})^2}$$

$$\leq \frac{1}{(k\delta^2/2)e^{-2\gamma t/L^2}} = \frac{2e^{2\gamma t/L^2}}{\delta^2 k}.$$

It follows that for $t \leq L^2 \log k / 4\gamma$ we have $\mathbf{P}(N_t \geq k/2) \leq 2/\delta^2 \sqrt{k}$, which is less than $\frac{1}{4}$ if $k > 64/\delta^4$. But, a uniformly chosen configuration will have at least $k/2$ black balls in $U$ with probability at least $\frac{1}{2}$. This completes the proof.

**Acknowledgments.** We are grateful to A. Bandyopadhyay, N. Berger, M. Liu, R. Lyons, F. Martinelli, E. Mossel, Y. Peres, A. Sinclair and P. Tetali for useful discussions.

MATHEMATICS DEPARTMENT
UNIVERSITY OF CALIFORNIA, DAVIS
1 SHIELDS AVENUE
DAVIS, CALIFORNIA 95616
USA
E-MAIL: morris@math.ucdavis.edu